% Logic Eprints
%Submitted 1059 Tue Oct 05, 1993 by: andreas.blass@math.lsa.umich.edu (andreas blass )
%logic/blass/divquot.tex
%

\magnification=1200
\input amstex
\documentstyle{amsppt}

%macros go here
\define\dv#1{\text{Div}(#1)}
\define\cf#1{\text{cf}(#1)}
\define\nn{{\Bbb N}}

\define\kw#1{{[#1]^\omega}}
\define\cp#1#2{{\text{CP}(#1\to#2)}}

\topmatter
\title
On the Divisible Parts of Quotient Groups
\endtitle
\author
Andreas Blass
\endauthor
\address
Mathematics Dept., University of Michigan, Ann Arbor, 
MI 48109, U.S.A.
\endaddress
\email
ablass\@umich.edu
\endemail
\thanks Partially supported by NSF grant DMS-9204276 and NATO grant
LG921395. 
\endthanks
\thanks
This paper is in final form, and no version of it will be submitted
for publication elsewhere.
\endthanks
\abstract
We study the possible cardinalities of the divisible part of $G/K$
when the cardinality of $K$ is known and when, for all countable
subgroups $C$ of $G$, the divisible part of $G/C$ is countable.
\endabstract
\endtopmatter
\document

\head
Introduction
\endhead

Let $G$ be an abelian group that is reduced, i.e., its divisible part
$\dv G$ is zero.  A quotient group $G/K$ need not be reduced.  This
paper is about the question how big $\dv{G/K}$ can be, relative to the
cardinality $|K|$ of $K$.

John Irwin has suggested, as a natural weakening of freeness, the
concept of a {\sl fully starred\/} torsion-free group, i.e., a
torsion-free, abelian group $G$ such that, for all subgroups $K$,
$|\dv{G/K}|\leq|K|$.  He asked whether, to test whether $G$ is fully
starred, it suffices to check the definition for countable $K$.  At
first sight, this seems unlikely; how should the quotients by
countable subgroups influence the quotients by larger subgroups.  We
shall show, however, that the answer to Irwin's question is
affirmative if $|K|<\aleph_\omega$.  We shall also show that it is
consistent with the usual axioms of set theory (ZFC, i.e.,
Zermelo-Fraenkel set theory, including the axiom of choice) that the
answer is negative for $|K|=\aleph_\omega$ but positive for many
larger values of $|K|$.  If a certain very large cardinal hypothesis
is consistent, then it is also consistent that the
answer to Irwin's question is affirmative for $|K|=\aleph_\omega$.  It
remains an open problem whether an affirmative answer, for all groups
regardless of cardinality, is consistent (relative to some large
cardinals).  

The following table summarizes the results.  The middle column gives
what can be said about $|\dv{G/K}|$ when $|K|$ is as in the first column,
when $\dv{G/C}$ is countable for all countable $C\leq G$, and when the
set-theoretic hypothesis in the right column is satisfied.  (The
assumption about countable $C$ is irrelevant in the first, very
elementary line of the table.)
\smallskip
$$\text{
\vbox{\offinterlineskip\hrule
\halign{
&\vrule#&\strut\hfil\ #\ \hfil\cr
height2pt&\omit&&\omit&&\omit&\cr
&$|K|=\kappa$&&$|\dv{G/K}|$&&Hypothesis&\cr
height2pt&\omit&&\omit&&\omit&\cr
\noalign{\hrule}
height2pt&\omit&&\omit&&\omit&\cr
&Arbitrary&&$\leq\kappa^{\aleph_0}$&&None&\cr
&$\kappa<\aleph_\omega$&&$\leq\kappa$&&None&\cr
&$\cf\kappa\neq\omega$&&$\leq\kappa$&&No inner model with measurable
cardinal&\cr
&$\cf\kappa=\omega$&&$\leq\kappa^+$&&No inner model with measurable
cardinal&\cr
&$\aleph_\omega$&&$\leq\aleph_\omega$&&Chang's conjecture for
$(\aleph_{\omega+1},\aleph_\omega)$&\cr
&$\aleph_\omega$&&possibly $=\aleph_{\omega+1}$&&$V=L$&\cr
&$\aleph_\omega$&&$<\aleph_{\aleph_4}$&&None&\cr
height2pt&\omit&&\omit&&\omit&\cr}
\hrule}}$$
\smallskip

Some of the hypotheses in the right column can be
weakened.  See Theorems~10 and 12 for sharper statements.

\head
1. Background
\endhead

Throughout this paper, all groups are abelian and reduced.  

In this section, we present, for motivation and orientation, some
elementary results and examples.  No novelty is claimed for any of
this material.  

As Irwin pointed out when he suggested the study of fully starred
torsion-free groups, this class of groups includes all free groups.
Indeed, if $G$ is freely generated by a basis $B$ and if $K$ is a
subgroup of $G$, then by expressing each element of $K$ as a
combination of (finitely many) elements of $B$, we obtain a subset
$B_0$ of $B$, no larger than $K$ in cardinality, such that $K$ is
included in the subgroup $G_0$ of $G$ generated by $B_0$.  Then
$G/K\cong (G_0/K)\oplus F$ where $F$ is freely generated by $B-B_0$.
Therefore the divisible part of $G/K$ coincides with that of $G_0/K$,
whose cardinality is at most that of $G_0$, which equals that of $K$.

For non-free $G$, on the other hand, the divisible part of $G/K$ may
well be larger than $K$.  For example, if $G$ is the product of
countably many infinite cyclic groups and $K$ is their direct sum
(embedded in the product in the obvious way), then the divisible part
of $G/K$ is easily seen to have the cardinality of the continuum, even
though $K$ is countable.  (Notice that, if the continuum hypothesis
holds, then the group $G$ in this example is almost free in the sense
that all subgroups of smaller cardinality are free \cite{10}.) Another
example is obtained by taking $G$ to be the additive group of $p$-adic
integers and $K$ the subgroup of ordinary integers; then $G/K$ has the
cardinality of the continuum and is divisible.

For countable $K$, the preceding examples achieve the largest possible
cardinality for $\dv{G/K}$, as the following proposition shows.

\proclaim{Proposition 1}
For any (reduced) group $G$ and any subgroup $K$, we have that
$|\dv{G/K}|\leq|K|^{\aleph_0}$. 
\endproclaim

\demo{Proof}
Let $D$ be the pre-image of $\dv{G/K}$ in $G$, so $K\subseteq
D\subseteq G$ and $D/K$ is divisible.  We shall prove the proposition
by producing a one-to-one function from $D$ into the set of countable
sequences of elements of $K$.  
As $D/K$ is divisible, we have, for
each $d\in D$ and each positive integer $n$, some $a_n(d)\in D$ and
some $b_n(d)\in K$ such that $d=n\cdot a_n(d)+b_n(d)$.  If $G$ is
torsion-free, then the map we
seek sends each $d\in D$ to the sequence $(b_n(d))_{n\in\nn}$.  To see
that it is one-to-one, suppose $d$ and $d'$ gave rise to the same
sequence.  Then, for each $n$, we would have
$d-d'=n\cdot(a_n(d)-a_n(d'))$. Thus, $d-d'$ belongs to the divisible
part of $G$ (here we use that $G$ is torsion-free), which is zero as
$G$ is reduced.

To avoid the assumption that $G$ is torsion-free, it suffices to
include, in the sequence associated to $d$, not only
all the $b_n(d)$ but also all $b_n(a_k(d))$, all $b_n(a_k(a_l(d)))$,
etc. Now, if $d$ and $d'$ give rise to the same sequence, then the
subgroup of $G$ generated by $d-d'$ and all elements of the forms
$a_k(d)-a_k(d')$, $a_k(a_l(d))-a_k(a_l(d'))$, etc. is divisible and
therefore zero.
\qed\enddemo

The proposition justifies the first line of the table in the
introduction.  As indicated there, this result, unlike the later ones,
does not depend on any assumption about $\dv{G/C}$ for countable $C$,
but only on the assumption (obviously needed) that $G$ is reduced.

We note that, by results of Cohen \cite2 and Solovay \cite9, the upper
bound given by Proposition~1 is rather weak.  Even when $K$ is
countable, the bound $|K|^{\aleph_0}$, which then equals the
cardinality of the continuum, can be an arbitrarily large cardinal.

\head
2. Some Combinatorial Set Theory
\endhead

We follow the usual set-theoretic conventions whereby an ordinal
number is the set of all smaller ordinal numbers and a cardinal number
is the first ordinal of that cardinality.  In particular, the cardinal
$\aleph_0$ is identified with the first infinite ordinal number
$\omega$ and with the set of all natural numbers.

For any set $A$, let $\kw A$ be the set of all countably infinite
subsets of $A$, partially ordered by the subset relation.  We shall
need some information about the cofinality $\cf{\kw A}$, i.e., the
smallest possible size for a family of countable subsets of $A$ such
that every countable subset of $A$ is included in one from that
family.  Of course this depends only on the cardinality of $A$, so we
may assume that $A$ is a cardinal.  The following proposition
summarizes the information we need about $\cf{\kw\kappa}$ for
uncountable $\kappa$; it is taken from Section 4 of \cite6.
The hypothesis ``the covering lemma over a model of GCH'' that occurs
in the second part of the proposition means that there is a class $M$
of sets such that every member of a member of $M$ is in $M$ (one says
$M$ is {\sl transitive\/}), all axioms of ZFC and the generalized
continuum hypothesis (GCH) hold in $M$, and every uncountable set of
ordinal numbers is included in a set of the same cardinality that is a
member of $M$.  We shall never need to use this definition of the
covering lemma over a model of GCH, but it is relevant that this
assumption is satisfied unless there are inner models with measurable
cardinals \cite{3,4}, i.e., transitive classes $M$ such that all axioms of
ZFC and the statement ``a measurable cardinal exists'' hold in $M$.
In fact, it is known that if this assumption fails then there are
inner models with far stronger large cardinal properties than just a
measurable cardinal; see for example \cite8.

\proclaim{Proposition 2}
The equation $\cf{\kw\kappa}=\kappa$ holds for all uncountable
cardinals $\kappa<\aleph_\omega$.  If the covering lemma holds over a
model of GCH, then the same equation holds for all cardinals $\kappa$
of uncountable cofinality, while for uncountable $\kappa$ of countable
cofinality $\cf{\kw\kappa}=\kappa^+$.\qed
\endproclaim

\demo{Proof}
See \cite6, Corollary~4.8, Lemma~4.10, and the proof of the latter.
\qed\enddemo

We note that the result obtained for the case of countable cofinality
using the covering lemma is optimal; that is, for $\cf\kappa=\omega$,
we cannot have $\cf{\kw\kappa}=\kappa$.  Indeed, given a family $\Cal
F$ of $\kappa$ countable subsets of $\kappa$ and given an increasing
$\omega$-sequence of infinite cardinals $\theta_i$ with supremum
$\kappa$, we can split $\Cal F$ into subfamilies $\Cal F_i$ of
respective cardinalities $\theta_i$, and we can find for each $i$ an
element $x_i\in\theta_{i+1}$ not in the union of $\Cal F_i$ (for this
union has size at most $\theta_i\cdot\aleph_0$).  Then the countable
set consisting of these $x_i$'s is not included in any member of $\Cal
F$.

In the absence of the covering lemma, it is much more difficult to
bound $\cf{\kw\kappa}$ for $\kappa$ of cofinality $\omega$.  The
following proposition is one of the surprising results of Shelah's
recently developed theory of possible cofinalities.  It is stated
(with a hint about the proof) as Proposition 7.13 of \cite1,
attributed to Baumgartner.

\proclaim{Proposition 3}
$\cf{\kw{\aleph_\omega}}<\aleph_{\aleph_4}$.
\endproclaim

The remainder of this section is devoted to a combinatorial principle
closely related to the preceding cofinality considerations but, as we
shall see, of more direct relevance to the size of the divisible parts
of quotient groups.

For infinite cardinals $\kappa<\lambda$, we say that the {\sl
$\lambda$ to $\kappa$ compression principle\/} holds and we write
$\cp\lambda\kappa$ if, for any $\lambda$-indexed family
$(X_i)_{i\in\lambda}$ of countable subsets of $\kappa$, there is an
uncountable set of indices, $Z\subseteq\lambda$, such that
$\bigcup_{i\in Z}X_i$ is countable.  Intuitively, this means that, if
$\lambda$ countable sets are packed (with overlapping) into a set of
size $\kappa$ (non-trivial packing, as $\kappa<\lambda$), then some
uncountably many of those sets must have been packed into a countable
set (so we have non-trivial packing on a smaller scale).

There is a trivial connection between the compression principle and
the cofinalities considered above; we state it as a proposition for
future reference.

\proclaim{Proposition 4}
If $\lambda>\cf{\kw\kappa}$ then $\cp\lambda\kappa$ holds.
\endproclaim

\demo{Proof}
Fix a family $\Cal F$ of strictly fewer than $\lambda$ countable
subsets of $\kappa$ cofinal in $\kw\kappa$, and fix $\lambda$
countable subsets $X_i$ of $\kappa$.  Each $X_i$ is in some member of
$\Cal F$.  As there are $\lambda$ $X_i$'s and fewer members of $\Cal
F$, uncountably many $X_i$'s (in fact $\lambda$ of
them) must be included in a single element of $\Cal F$ and must
therefore have a countable union.
\qed\enddemo

The converse of Proposition~4 is not provable, at least if
sufficiently large cardinals are consistent.  The counterexample
involves Chang's conjecture.  Originally, Chang's conjecture was that,
for a countable first-order language,
every structure of cardinality $\aleph_2$ in which a
particular unary predicate symbol $P$ denotes a set of cardinality
$\aleph_1$ has an elementary submodel of cardinality $\aleph_1$ in
which $P$ denotes a set of cardinality $\aleph_0$.  Generalizing this
by changing $\aleph_2$ and $\aleph_1$ in the hypothesis to $\lambda$
and $\kappa$, respectively, but leaving the conclusion unchanged, one
has the Chang conjecture for $(\lambda,\kappa)$.  The conjecture can
be restated in a form, more convenient for our purposes, that avoids
model-theoretic notions.  We adopt this restatement as the definition:
{\sl Chang's conjecture\/} for $(\lambda,\kappa)$ is the assertion
that, for any countably many functions $f_n$, each mapping some finite
power $\lambda^p$ of $\lambda$ into $\kappa$ (where $p$ can depend on
$n$), there is an uncountable subset $H$ of $\lambda$ such that each
$f_n$ restricted to $H$ (more precisely, to $H^p$) has countable
range.  We shall need an elementary connection between Chang's
conjecture and the compression principle and a (non-elementary)
theorem from \cite7 giving the consistency of a particular instance of
Chang's conjecture.  

\proclaim{Proposition 5}
Chang's conjecture for $(\lambda,\kappa)$ implies $\cp\lambda\kappa$. 
\endproclaim

\demo{Proof}
Assume Chang's conjecture for $(\lambda,\kappa)$, and let $\lambda$
countable subsets $X_i,\ i<\lambda$ of $\kappa$ be given.  Fix, for
each $i$ an enumeration of $X_i$ by natural numbers, and define
functions $f_n:\lambda\to\kappa$ for $n\in\omega$ by letting $f_n(i)$
be the $n$th element in the chosen enumeration of $X_i$.  By Chang's
conjecture for $(\lambda,\kappa)$, there is an uncountable
$Z\subseteq\lambda$ whose images under all the $f_n$ are countable.
Then $\bigcup_{i\in Z}X_i$, which is also the union of these countably
many countable images, is countable, as required by the compression
principle.
\qed\enddemo

To state concisely the consistency theorem for Chang's conjecture for
$(\lambda,\kappa)=(\aleph_{\omega+1},\aleph_\omega)$, we introduce a
short name for the large cardinal hypothesis needed.  (We shall not
explicitly use this hypothesis, so the reader can safely skip the
definition and remember only that we are dealing with very large
cardinals, much larger than measurable cardinals but not known or even
widely believed to be inconsistent.)  A cardinal $\kappa$ will be
called {\sl huge\/}+ if there is an elementary embedding $j$ of the
set-theoretic universe into a transitive class $M$ such that $\kappa$
is the first ordinal moved and, if $\mu$ denotes the $(\omega+1)$th
cardinal after $j(\kappa)$ then every function from $\mu$ into $M$ is
an element of $M$.  ({\sl Huge\/} is defined the same way except that
$\mu=j(\kappa)$.)

\proclaim{Proposition 6}
If the existence of a huge{\rm+} cardinal is consistent with ZFC, then so is
Chang's conjecture for $(\aleph_{\omega+1},\aleph_\omega)$.
\endproclaim

\demo{Proof}
See Levinski, Magidor, and Shelah \cite7, Theorem 5.
\qed\enddemo

We conclude this section with an analog for the compression property
of a well-known fact about $\cf{\kw\kappa}$ \cite6, Lemma 4.6.

\proclaim{Proposition 7}
If $\lambda>\kappa^+$ and if $\cp\lambda\kappa$ then
$\cp\lambda{\kappa^+}$.  More generally, if $\lambda>\mu$, if
$\cf\mu>\omega$, and if $\cp\lambda\kappa$ for all $\kappa<\mu$, then
$\cp\lambda\mu$.  
\endproclaim

\demo{Proof}
Since $\cp\lambda\kappa$ trivially implies $\cp\lambda\theta$ for all
$\theta<\kappa$, it suffices to prove the second statement.  Let
$\lambda$ countable subsets $X_i$ of $\mu$ be given.  As
$\cf\mu>\omega$, each $X_i$ has its supremum $<\mu$, and, as
$\lambda>\mu$, this supremum must be the same ordinal, say $\alpha$,
for $\lambda$ of the $X_i$'s.  But then $\cp\lambda{|\alpha|}$ ensures
that some uncountably many of these $X_i$ have a countable union.
\qed\enddemo

\head
3. Divisible Parts of Quotients Are Not Too Large
\endhead

In this section, we present our positive results about Irwin's
question described in the introduction.  That is, we assume that
$|\dv{G/C}|$ is countable for all countable subgroups $C$ of $G$, and
we deduce upper bounds for $|\dv{G/K}|$ in terms of $|K|$.  All these
results are obtained by combining the set-theoretic results in
Section~2 with the following proposition which relates the divisible
parts of quotient groups to the compression principle.

\proclaim{Proposition 8}
Let $G$ be a group such that $\dv{G/C}$ is countable for all countable
subgroups $C$ of $G$.  Let $\lambda>\kappa$ be cardinals such that
$\cp\lambda\kappa$ holds.  Then for all subgroups $K$ of $G$ of
cardinality $\kappa$, $|\dv{G/K}|<\lambda$.
\endproclaim

\demo{Proof}
Assume toward a contradiction that the hypotheses hold but the
conclusion fails.  Then there is a subgroup $D$ of cardinality
$\lambda$ with $K\subseteq D\subseteq G$ and with $D/K$ divisible.  As
in the proof of Proposition~1, associate to each $d\in D$ and each
positive integer $n$ elements $a_n(d)\in D$ and $b_n(d)\in K$ such
that $d=n\cdot a_n(d)+b_n(d)$.  For each $d\in D$, let $X_d$ be the
set of those elements of $K$ obtainable by applying to $d$ any finite
composite of the functions $a_n:D\to D$ followed by any of the
functions $b_n:D\to K$.  (The identity function counts as a finite
composite of $a_n$'s, namely the composite of none.)  Thus $X_d$
consists of all elements of the forms $b_n(d)$, $b_n(a_k(d))$,
$b_n(a_k(a_l(d)))$, etc.  

As each $X_d$ is a countable subset of the $\kappa$-element set $K$,
the assumed $\lambda$ to $\kappa$ compression principle provides an
uncountable subset $Z$ of $D$ such that all $X_d$ for $d\in Z$ are
contained in a countable subset $C$ of $K$.  We may assume that $C$ is
a subgroup of $K$, by replacing it with the subgroup it generates.  

Let $E$ be the subgroup of $D$ generated by all the elements
obtainable from elements of $Z$ by applying any finite composite of
the functions $a_n$.  Notice that, if $e$ is any one of these
generators, then $b_n(e)\in C$ and $a_n(e)\in E$ for all positive
integers $n$.  Thus, from $e=n\cdot a_n(e)+b_n(e)$, we infer that
$E/C$ is divisible.  But, as $E$ is uncountable (containing $Z$) and
$C$ is countable, $E/C$ is an uncountable subgroup of $G/C$, contrary
to the assumption that, for countable $C$, $\dv{G/C}$ is countable. 
\qed\enddemo

\proclaim{Corollary 9}
Let $G$ be a group such that $\dv{G/C}$ is countable for all countable
subgroups $C$ of $G$.  If $K$ is a subgroup of $G$ of cardinality
$\kappa$, then $|\dv{G/K}|\leq\cf{\kw\kappa}$.
\endproclaim

\demo{Proof}
Combine Propositions~4 and 8.
\qed\enddemo

\proclaim{Theorem 10}
Let $G$ be a group such that $\dv{G/C}$ is countable for all countable
subgroups $C$ of $G$, and let $K$ be a subgroup of $G$ of cardinality
$\kappa$. 
\roster
\item If $\kappa<\aleph_\omega$ then $|\dv{G/K}|\leq\kappa$.
\item If the covering lemma holds over a model of GCH (in particular
if there is no inner model with a measurable cardinal) and if
$\cf\kappa>\omega$, then $|\dv{G/K}|\leq\kappa$.
\item If the covering lemma holds over a model of GCH (in particular
if there is no inner model with a measurable cardinal) and if
$\cf\kappa=\omega$, then $|\dv{G/K}|\leq\kappa^+$.
\item If Chang's conjecture for $(\aleph_{\omega+1},\aleph_\omega)$ is
true and if $\aleph_\omega\leq\kappa<\aleph_{\omega\cdot2}$, then 
$|\dv{G/K}|\leq\kappa$.
\item If $\kappa=\aleph_\omega$ then $|\dv{G/K}|<\aleph_{\aleph_4}$.
\endroster
\endproclaim

\demo{Proof}
\therosteritem1, \therosteritem2, and \therosteritem3 all follow
immediately from Proposition~2 and Corollary~9. For \therosteritem4,
use Propositions~5 and 7 to deduce from Chang's conjecture for
$(\aleph_{\omega+1},\aleph_\omega)$ that $\cp\kappa^+\kappa$ holds for
all $\kappa$ as in \therosteritem4; then invoke Proposition~8.
\therosteritem5 follows from Proposition~3 and Corollary~9.
\qed\enddemo

It follows from \therosteritem1, \therosteritem4 and Proposition~6
that the statement ``If $G$ is a group such that $\dv{G/C}$ is
countable for all countable subgroups $C$ of $G$, and if $K$ is a
subgroup of $G$ with $|K|<\aleph_{\omega\cdot2}$, then
$|\dv{G/K}|\leq|K|$'' is consistent relative to a huge+
cardinal.

The results presented in this section complete the justification of
all lines in the table in the introduction except for the next to last
line.  The one remaining line, asserting that the answer to Irwin's
question is negative for groups of cardinality $\aleph_\omega$ if
$V=L$, is the subject of the next two sections.

\head
4. More Combinatorial Set Theory
\endhead

In this section, we develop, under the assumption $V=L$ the set theory
needed to produce a counterexample for Irwin's question.  In fact, we
do not need the full strength of $V=L$ but only the combinatorial
principle $\square_{\aleph_\omega}$.  For any uncountable cardinal
$\kappa$, $\square_\kappa$ denotes the following assertion: There
exists a sequence of sets $(C_\xi)$ indexed by the limit ordinals
$\xi<\kappa^+$ such that for each such $\xi$
\roster
\item $C_\xi$ is a closed, cofinal subset of $\xi$.
\item If $\cf\xi<\kappa$ then $|C_\xi|<\kappa$.
\item If $\eta$ is a limit point of $C_\xi$, then $C_\eta=\eta\cap
C_\xi$. 
\endroster
These square principles were introduced by Jensen \cite5 who showed
that they follow from $V=L$ for all $\kappa$.  

In order to obtain a negative answer to Irwin's question at
cardinality $\aleph_\omega$, it is necessary, according to
Proposition~8, to violate $\cp{\aleph_{\omega+1}}{\aleph_\omega}$,
i.e., to produce $\aleph_{\omega+1}$ countable subsets of
$\aleph_\omega$ such that no countable set contains uncountably many
of them.  The following proposition shows that this and a bit more
(which we shall need in the next section) can be done if
$\square_{\aleph_\omega}$ holds.  I am not sure to whom to attribute
this proposition.  Menachem Kojman told me that
$square_{aleph_\omega}$ contradicts
$\cp{\aleph_{\omega+1}}{\aleph_\omega}$.  Menachem Magidor showed me
the proof given below for this fact and the additional information in
Proposition~11.  Magidor also informed me that similar arguments were
known to Saharon Shelah.

\proclaim{Proposition 11}
Assume $\square_{\aleph_\omega}$.  There are $\aleph_{\omega+1}$
countably infinite subsets $X_\xi$ of $\aleph_\omega$ such that the
intersection of every two of these $X_i$ is finite and, for
each countable $Y\subseteq\aleph_\omega$, at most countably many of the
$X_\xi$ have infinite intersection with $Y$ (and, a fortiori, at most
countably many $X_\xi$ are included in $Y$).
\endproclaim

\demo{Proof}
Fix $C_\xi$ for limit ordinals $\xi<\aleph_{\omega+1}$ as in the
definition of $\square_{\aleph_\omega}$.  We shall define functions
$f_\xi$ for $\xi<\aleph_{\omega+1}$ with the following properties.
\roster
\item Each $f_\xi$ is a function on $\omega$ satisfying
$f_\xi(n)\in\aleph_n$ for all $n\in\omega$.
\item If $\xi<\eta$, then $f_\xi(n)<f_\eta(n)$ for all but finitely many
$n\in\omega$. 
\item If $\eta$ is a limit ordinal, $|C_\eta|<\aleph_n$, and $\xi\in 
C_\eta$, then $f_\xi(n)<f_\eta(n)$.
\endroster
(Notice that, for $\xi\in C_\eta$, \therosteritem3 amplifies
\therosteritem2 by specifying that the finitely many exceptional $n$
in \therosteritem2 are bounded above by the $q$ such that
$|C_\eta|=\aleph_q$.  Notice also that such a $q$ exists by
\therosteritem3 in the definition of $\square_{\aleph_\omega}$.)  After
constructing such $f_\xi$'s, we shall show that their graphs are
essentially the sets needed to establish the proposition.

The construction of the $f_\xi$'s proceeds by induction on
$\xi<\aleph_{\omega+1}$.  Suppose, therefore, that $f_\xi$ is defined
for every $\xi<\eta$, and we wish to define $f_\eta$.

If $\eta$ is not a limit ordinal, then
\therosteritem3 does not apply to $\eta$, so we need only satisfy
\therosteritem1 and \therosteritem2, which we do as follows.
Partition $\eta$ (the set of ordinals smaller than $\eta$), which has
cardinality at most $\aleph_\omega$, into countably many pieces $A_0,
A_1, \dots$ such that each $A_k$ has cardinality at most $\aleph_k$.
Then define $f_\eta(n)$ to be any ordinal that is $<\aleph_n$ (so
\therosteritem1 holds) but $>f_\xi(n)$ for all $\xi\in
A_0\cup\dots\cup A_{n-1}$ (so that \therosteritem2 holds, because if
$\xi\in A_p$ then $f_\xi(n)<f_\eta(n)$ for all $n>p$).  Such an
ordinal exists, because there are only $\aleph_{n-1}$ ordinals $\xi\in
A_0\cup\dots\cup A_{n-1}$ and therefore the corresponding $f_\xi(n)$'s
cannot be cofinal in $\aleph_n$.  

If $\eta$ is a limit ordinal, let $q$ be the natural number such that
$|C_\eta|=\aleph_q$.  Then, in order to satisfy \therosteritem3, we
must make sure that $f_\eta(n)$ satisfies, in addition to the
requirements in the preceding paragraph, $f_\xi(n)<f_\eta(n)$ if
$\xi\in C_\eta$ and $n>q$.  But there are only $\aleph_q$ such
$\xi$'s, by choice of $q$, so the corresponding $f_\xi(n)$'s cannot be
cofinal in $\aleph_n$ when $n>q$.  So an appropriate value for
$f_\eta(n)$ can again be found.  This completes the construction of
the $f_\eta$'s and the verification of \therosteritem1,
\therosteritem2, and \therosteritem3.

We claim that, from any uncountable $A\subseteq\aleph_{\omega+1}$, one
can extract an uncountable $B\subseteq A$ and one can find
$m\in\omega$ such that $f_\xi(n)<f_\eta(n)$ for all $\xi<\eta$ in $B$
and all $n\geq m$.  (That is, the finitely many exceptions in
\therosteritem2 are bounded by the same $m$, independent of $\xi$ and
$\eta$, as long as these two indices lie in $B$.)  To prove this
claim, consider an arbitrary uncountable
$A\subseteq\aleph_{\omega+1}$.  We may assume without loss of
generality that $A$ has order type $\aleph_1$; just replace $A$ by the
subset consisting of its first $\aleph_1$ elements.  Let $\beta$ be
the supremum of $A$, so $\cf\beta=\aleph_1$.  Inductively choose
$$
\gamma(0)<\alpha(0)<\gamma(1)<\alpha(1)<\dots
<\gamma(\omega)<\alpha(\omega)<\gamma(\omega+1)<\dots
$$
for $\aleph_1$ steps, so that all the $\gamma$'s are limit points of
$C_\beta$ and all the $\alpha$'s are in $A$.  There is no difficulty
making these choices, as both the set of limit points of $C_\beta$ and
$A$ are cofinal in $\beta$.  For each $\mu<\aleph_1$,
\therosteritem2 provides a natural number $m(\mu)$ such that
$$
\forall n\geq m(\mu)
\quad f_{\gamma(\mu)}(n)<f_{\alpha(\mu)}(n)<f_{\gamma(\mu+1)}(n).
$$
Increasing $m(\mu)$ if necessary, we may assume that
$|C_\beta|<\aleph_{m(\mu)}$.  
As there are uncountably many $\mu$'s and only countably many
possible values for $m(\mu)$, there is an uncountable
$Y\subseteq\aleph_1$ such that $m(\mu)$ has the same value $m$ for
all $\mu\in Y$.  Now if $\mu<\nu$ are in $Y$ and if $n\geq m$ then 
we have
$$
f_{\alpha(\mu)}(n)<f_{\gamma(\mu+1)}(n)\leq
f_{\gamma(\nu)}(n)<f_{\alpha(\nu)}(n),
$$
where the first and third inequalities hold because $n\geq
m=m(\mu)=m(\nu)$, and the middle inequality holds because of
\therosteritem3 and the fact that, since the $\gamma$'s are limit
points of $C_\beta$, either $\mu+1=\nu$ or
$\gamma(\mu+1)\in\gamma(\nu)\cap C_\beta=C_(\nu)$. This 
means that $B=\{\alpha(\mu)\mid\mu\in Y\}$ and $m$ are as required in
the claim.

Finally, we show that the functions $f_\xi$ (for all
$\xi\in\aleph_{\omega+1}$), regarded as sets of ordered pairs
$\subseteq\omega\times\aleph_\omega$, have the properties that every
two have finite intersection and that no
countable set $C$ can have infinite intersections with uncountably
many $f_\xi$'s.  Then, transferring these subsets of
$\omega\times\aleph_\omega$ to subsets of $\aleph_\omega$ by some
bijection, we have the sets required by the proposition.  That every
two of the graphs have finite intersection is immediate from
\therosteritem2.  To prove the other property, suppose
$C$ were a countable set having infinite intersection with $f_\xi$ for
all $\xi$ in some uncountable $A\subseteq\aleph_{\omega+1}$.  Let
$B\subseteq A$ and $m$ be as in the claim proved in the preceding
paragraph.  Then by deleting the first $m$ elements from each of the
graphs $f_\xi$, $\xi\in B$, i.e., by forming
$f_\xi\restriction(\omega-m)$, we would obtain uncountably many pairwise
disjoint sets, all intersecting the countable set $C$.  As this is
absurd, the proof is complete.
\qed\enddemo

\head
5. A Quotient with a Large Divisible Part
\endhead

In this section we construct, assuming $\square_{\aleph_\omega}$, a
counterexample to Irwin's question, i.e., a group $G$ such that
$\dv{G/C}$ is countable for all countable $C$ but $G$ is not fully
starred.  

\proclaim{Theorem 12}
Assume $\square_{\aleph_\omega}$.  There exist a group $G$ of
cardinality $\aleph_{\omega+1}$ and a subgroup $K$ of cardinality
$\aleph_\omega$ such that $G/K$ is divisible but, for each countable
subgroup $C$ of $G$, the divisible part of $G/C$ is countable.
\endproclaim

\demo{Proof}
The required groups $G$ and $K$ will actually be modules over the ring
$R$ of rational numbers with odd denominators, i.e., they will be
divisible by all primes except 2.  We write $R^*$ for the set of units
of $R$, the rational numbers whose numerators and denominators (in
reduced form) are both odd.

Since $\square_{\aleph_\omega}$ is assumed, let $X_i$ for
$i\in\aleph_{\omega+1}$ be as in Proposition~11.  Fix, for each $i$ a
bijection between $X_i$ and the set $\Bbb Z$ of integers, and write
$\xi(i,n)$ for the element of $X_i$ that corresponds to the integer
$n$.  (Notice that the same ordinal can be $\xi(i,n)$ for many
different pairs $(i,n)$.)

The group $G$ will be presented as the $R$-module generated by certain
objects subject to certain relations.  The generators are of two
sorts.  First, there are $\aleph_\omega$ generators which, to simplify
notation, we take to be the ordinals $\alpha<\aleph_\omega$.  Second,
there are $\aleph_{\omega+1}$ generators $g(i,n)$ indexed by all the
ordinals $i<\aleph_{\omega+1}$ and all integers $n\in\Bbb Z$.  The
defining relations are 
$$
g(i,n)=2\cdot g(i,n+1)+\xi(i,n).
$$
This defines $G$ as an $R$-module, hence as a group.  $K$ is defined
to be the submodule generated by the first sort of generators of $G$,
the ordinals below $\aleph_\omega$.

Notice that, in any non-trivial $R$-linear combination of the defining
relations of $G$, some $g(i,n)$ must occur.  Indeed, of the finitely
many $g(i,n)$'s that occur in the relations being combined, the ones
with the largest (or the smallest) values of $n$ cannot be canceled.
So the relations impose no restrictions on the generators of $K$
alone.  This means that $K$ is freely generated, as an $R$-module, by
the ordinals $\alpha<\aleph_\omega$.  It follows, exactly as in the
argument for free groups in Section~1, that $K$ is fully starred.

The quotient group $G/K$ can be presented by adjoining to the defining
relations of $G$ the new relations $\alpha=0$ for all the generators
$\alpha$ of $K$.  The resulting presentation is clearly equivalent to
one that has only the generators $\bar g(i,n)$ (where the bar over a
letter denotes the coset modulo $K$) and the relations
$$
\bar g(i,n)=2\cdot \bar g(i,n+1).
$$  
Thus, $G/K$ is divisible and is in fact
the rational vector space freely generated by the $\bar g(i,0)$'s; here
$\bar g(i,n)$ is identified with $2^{-n}\bar g(i,0)$.

To complete the proof of the theorem, it remains to show that
$\dv{G/C}$ is countable for all countable subgroups $C$ of $G$.  As a
preliminary step toward this, we introduce a normal form for elements
of $G$.  We claim that every element $x$ of $G$ can be uniquely
represented in the form
$$
x=\sum_ir_ig(i,n_i)+\sum_\alpha s_\alpha\alpha
$$
where all $r_i\in R^*$, $s_\alpha\in R$, and the summation variables
$i$ and $\alpha$ range over some finite
subsets (depending on $x$) of $\aleph_{\omega+1}$ and $\aleph_\omega$,
respectively.  (In particular,
each $i$ that occurs at all, as the first argument of a $g$,  occurs 
in only one term of such a normal form.)  To produce such a normal
form for $x$, we first write the image of $x$ in $G/K$ as a rational
linear combination of the $\bar g(i,0)$'s, then we put the rational
coefficients in this combination into $R^*$ by using $2^n
\bar g(i,0)=\bar g(i,-n)$, which holds in $G/K$, then we remove the
bars from the $g$'s in this expression, obtaining an element $y\in G$,
and we observe that the difference $d=x-y$ is in $K$ as $y$ has the
same image as $x$ in $G/K$.  Now the desired normal form of
$x$ consists of $y$ plus the expansion of $d$ in terms of the free
generators $\alpha$ of $K$.  The uniqueness is proved by observing that the
$\sum_ir_ig(i,n_i)$ part of a normal form of $x$ must have the same
image as $x$ in $G/K$ and must therefore be exactly the $y$
constructed above.

It will be useful later to have a more computational description of
how to convert an arbitrary element of $G$, given as an $R$-linear
combination of generators $g(i,n)$ and $\alpha$, to normal form.  This
means that we must convert the given expression, using the defining
relations of $G$, to one in which, for each $i$, there is at most one
term of the form $r\cdot g(i,n)$ and the coefficient $r$ of any such
term is in $R^*$.  Suppose that, for a certain $i$, the given
expression contains several terms of the form $g(i,n)$, possibly with
different values of $n$.  The defining relations of $G$ allow us to
replace $g(i,n)$ with $g(i,n+1)+\xi(i,n)$, so we can increase the
$n$'s involved in these terms as much as we wish.  In particular, we
can increase them until they are all equal to, say, the largest of the
originally occurring $n$'s.  At this stage, all the $g(i,n)$'s, for
the particular $i$ under consideration, have the same $n$, so we can
collect them into one term by adding their coefficients.  If the
resulting coefficient of $g(i,n)$ is divisible by 2, then we can use
the defining relations ``in reverse'' to divide the coefficient by 2
while decreasing $n$ by 1 and subtracting $\xi(i,n-1)$.  Repeat this
until the coefficient is no longer divisible by 2, i.e., until it is
in $R^*$.  By carrying out the procedure just described for each $i$,
we clearly achieve normal form.

We point out two consequences of this procedure that will be useful
later.  First, notice that the element $\sum_\alpha s_\alpha\alpha\in
K$  occurring in the normal
form at the end consists of first the linear combination of
$\alpha$'s in the original expression, and second some terms of the
form $\pm2^k\xi(i,n)$, where the $i$'s involved in these $\xi$'s were
involved in $g$'s in the original expression.

Second, if $x-2z\in K$ then $g(i,n)$ occurs in the normal form
of $x$ if and only if $g(i,n+1)$ occurs in the normal form of $z$, and
they occur with the same coefficient.  To see this, start with a
normal form of $z$, multiply all terms by 2 and add a suitable element
of $K$ to get an expression for $x$ that is in normal form except that
all the coefficients of $g$'s are divisible once by 2; then apply the
normalization algorithm to this expression.

We now embark on the proof that $\dv{G/C}$ is countable for all
countable $C\subseteq G$.  Notice first that, if $C$ were a
counterexample and if $C'$ were any countable group such that
$C\subseteq C'\subseteq G$, then $C'$ would also be a counterexample.
Indeed, there would be an uncountable $D$ with $C\subseteq D\subseteq
G$ and $D/C$ divisible.  Then, as $C\subseteq D\cap C'$,
$(D+C')/C'\cong D/(D\cap C')$ is a quotient of the divisible group
$D/C$, and is therefore divisible.  As $D+C'$ is uncountable and $C'$
is countable, we have another counterexample, as claimed.

Therefore, in proving that $G/C$ is divisible, we may assume that $C$
satisfies the following four conditions, since each condition 
amounts to being
closed under countably many functions and can therefore be satisfied
by a suitable countable supergroup of any given countable $C$.
\roster
\item $C$ is an $R$-submodule of $G$.
\item If $x\in C$, then all the generators $g(i,n)$ and $\alpha$ that
occur in the normal form of $x$ are also in $C$.
\item If $g(i,n)\in C$, then $g(i,m)\in C$ for all integers $m$.
\item If $g(i,n)\in C$, then all members $\xi(i,m)$ of $X_i$ are in
$C$.  
\endroster
We call $i<\aleph_{\omega+1}$ a $C$-{\sl index\/} if, for some $n$,
$g(i,n)\in C$.  By \therosteritem2, it is equivalent to say that
$g(i,n)$ occurs in the normal form of some element of $C$, and by
\therosteritem3 it is also equivalent to say that $g(i,n)\in C$ for
all $n\in \Bbb Z$.  As $C$ is countable, there are only countably many
$C$-indices.  (Note, by contrast, that there may be uncountably many
$i$ such that, for some $n$, $\xi(i,n)\in C$, for the same element of
$C$ can be $\xi(i,n)$ for many different $i$.)

Now suppose $D$ is a subgroup of $G$ with $C\subseteq D$ and with
$D/C$ divisible.  We must prove that $D$ is countable.  Call an
element $i\in\aleph_{\omega+1}$ {\sl relevant\/} if the normal form of
some element of $D$ contains $g(i,n)$ for some $n$.

\proclaim{Lemma}
Only countably many $i$ are relevant.
\endproclaim

\demo{Proof}
Suppose uncountably many $i$ are relevant.  Then the set $U$ of
relevant $i$ that are not $C$-indices is uncountable.  As the $X_i$
were chosen as in 
Proposition~11, the countable set $C\cap\aleph_\omega$ (i.e., the set
of generators of $K$ that are in $C$) must have finite intersection
with $X_i$ for some $i\in U$.  Fix such an $i$ for the rest of the
proof of the lemma.  Thus, $i$ is relevant, $i$ is not a $C$-index,
and $C\cap X_i$ is finite.

As $i$ is relevant, choose an element $d_0\in D$ such that $g(i,q)$
occurs in $d_0$ for some $q$.  (We think of $d_0$ and the other
elements defined below as written in normal form, so ``occurs in $d_0$''
really means ``occurs in the normal form of $d_0$.'')  Let $a\in R^*$
be the coefficient of $g(i,q)$ in $d_0$.  Let $J$ be the
set of those $j$ such that (i) $g(j,m)$ occurs in $d_0$ for some $m$,
(ii) $j\neq i$, and (iii) $j$ is not a $C$-index.  Of course, $J$ is
finite. 

Since $X_i$ has finite intersection with $C$ (by choice of $i$) and
with $X_j$ for each $j\in J$ (by the ``finite pairwise intersections''
part of Proposition~11), fix $p\in\Bbb Z$ so large that, for all
$n\geq p$, we have $\xi(i,n)\notin C$ and $\xi(i,n)\notin X_j$ for all
$j\in J$.  

Define a sequence $(d_r)_{r\in\omega}$ of elements of $D$ by starting
with $d_0$ and then inductively using the
divisibility of $D/C$ to set $d_r=2\cdot d_{r+1}+c_r$ with $d_{r+1}\in
D$ and $c_r\in C$.  

We claim that (the normal form of) $d_r$ contains $g(i,q+r)$ with
coefficient $a$. We chose $a$ to make this true for $r=0$.  If it is
true for $r$, then, as $g(i,n)$ does not occur in $c_r$ for any $n$
(as $i$ is not a $C$-index), $g(i,q+r)$ occurs with coefficient $a$ in
$d_r-c_r=2\cdot d_{r+1}$.  By the second
consequence of the algorithm for converting expressions to normal
form, it follows that $g(i,q+r+1)$ occurs with coefficient $a$ in the
normal form of $d_{r+1}$, as desired.

Again referring to the algorithm for computing normal forms, we see
that the set $J$, defined above using $d_0$, would be unchanged if we
used any $d_r$ instead of $d_0$ in clause (i).  Indeed, if some $j$
satisfied clause (i) for one of $d_r$
and $d_{r+1}$ but not for the other, then the equation $d_r=2\cdot
d_{r+1}+c_r$ would require some $g(j,m)$ to occur in $c_r$, so
$j$ would be a $C$-index and would thus fail to satisfy clause (iii).

In view of the observations in the preceding two paragraphs, any final
segment of the sequence $(d_r)$ is again a sequence of the same sort,
only with a larger value of $q$.  Since $p$ can obviously also be
replaced by any larger number, we assume without loss of generality
that $p=q$.  Thus, henceforth, $d_r$ has $g(i,p+r)$ in its normal
form, with coefficient $a\in R^*$.

The definition of the $d$ sequence gives us that 
$$\align
d_0&=2d_1+c_0\\
&=4d_2+2c_1+c_0\\
&=\dots\\
&=2^rd_r+2^{r-1}c_{r-1}+\dots+2c_1+c_0,
\endalign$$
for any $r$.
Fixing some $r\geq 1$, let us consider the normal form of $d_0$ as
obtained from the last line in the preceding display by first writing
$d_r$ and $c=2^{r-1}c_{r-1}+\dots+2c_1+c_0$ in normal form, then
combining these normal forms to make an expression for $d_0=2^rd_r+c$,
and finally normalizing the result.  More precisely, let us consider
the coefficient $s\in R$ of $\xi(i,p+r-1)$ in the resulting normal
form, and in fact let us consider $s$ modulo $2^r$.  

First, let us consider possible occurrences of $\xi(i,p+r-1)$ in the
first expression for $2^rd_r+c$ that we got by combining normal forms
of $d_r$ and $c$ (before normalizing the result).  Any ordinal
$\alpha$ that occurs in the normal form of $c\in C$ is itself in $C$
(closure condition \therosteritem2), and hence is different from
$\xi(i,p+r-1)$ by our choice of $p$.  There may be occurrences of
$\xi(i,p+r-1)$ in $d_r$, but these will have their coefficients
multiplied by $2^r$ in the expression for $2^rd_r+c$, and therefore
will not contribute to the coefficient $s$ modulo $2^r$.

It remains to consider occurrences of $\xi(i,p+r-1)$ that arise during
the normalization of $2^rd_r+c$.  These can arise as $\xi(j,n)$ during
the normalization of $g(j,m)$ terms.  (Recall that $\xi$'s with
different arguments can be equal, so we cannot a priori exclude cases
with $j\neq i$.  The rest of this paragraph shows that these cases
can, nevertheless, be excluded.)  Consider first those $g(j,m)$ for
which $j$ is a $C$-index.  That includes all the $g$'s in the normal
form of $c$ and possibly some in the normal form of $d_r$ as well.
Any $\xi(j,n)$ arising during the normalization of these $g$'s is in
$C$ by closure condition \therosteritem4 on $C$.  These $\xi(j,n)$'s
therefore differ from $\xi(i,p+r-1)$ by our choice of $p$.  It remains
to consider those $g(j,m)$ in the normal form of $d_r$ that are not
$C$-indices.  One of these has $j=i$, which we treat in the next
paragraph.  The rest have $j\in J$. Any $\xi$ produced by their
normalization is therefore in $X_j$ with $j\in J$ and is therefore
different from $\xi(i,p+r-1)$ by our choice of $p$.

We have seen that the coefficient of $s$ modulo $2^r$ must arise
entirely from the normalization of the term $2^rag(i,p+r)$ in
$2^rd_r+c$. This normalization process reads 
$$\align
2^rag(i,p&+r)=2^{r-1}ag(i,p+r-1)-2^{r-1}a\xi(i,p+r-1)\\
&=2^{r-2}ag(i,p+r-2)-2^{r-2}a\xi(i,p+r-2)-2^{r-1}a\xi(i,p+r-1)\\
&=\dots\\
&=ag(i,p)-a\xi(i,p)-2a\xi(i,p+1)-\dots-2^{r-1}a\xi(i,p+r-1).
\endalign$$
So we see that the coefficient of $\xi(i,p+r-1)$ is $2^{r-1}a$.  Since
$a\in R^*$, this coefficient is not zero modulo $2^r$.

This proves that the coefficient of $\xi(i,p+r-1)$ in the normal form
of $d_0$ is not divisible by $2^r$ and, in particular, is not zero.
But $r$ was an arbitrary integer $\geq1$, so the normal form of $d_0$
contains infinitely many non-zero terms.  That is absurd, and so the
lemma is proved.
\qed\enddemo

It follows from the lemma just proved that the image of the projection
of $D$ to $G/K$ is countable, for it is generated as an $R$-module by
$\bar g(i,n)$'s where $i$ is relevant and $n\in\Bbb Z$.  That image is
$D/(D\cap K)$, so to
complete the proof that $D$ is countable and thus the proof of the
theorem, it suffices to show that the kernel $D\cap K$ of that
projection is also countable.  We remarked earlier that $K$, being a
free $R$-module, is fully starred, so we need only prove that $(D\cap
K)/(C\cap K)$ is divisible.  This we now do.

Consider an arbitrary element $d$ of $D\cap K$.  As $D/C$ is
divisible, we have $d=2y-c$ for some $y\in D$ and some $c\in C$.  Let
their normal forms be
$$\align
c&=\sum_ir_ig(i,n_i)+k_c\\
y&=\sum_ir_ig(i,n_i+1)+k_y,
\endalign$$
where $k_c,k_y\in K$ and where the $g$ terms match except for a shift
of $n$'s because $c-2y\in K$. (See the second consequence of the
normalization algorithm.)  Thus,
$$\align
d=2y-c&=\sum_ir_i[2g(i,n_i+1)-g(i,n_i)]+2k_y-k_c\\
&=-\sum_i\xi(i,n_i)+2k_y-k_c.
\endalign$$
As $c\in C$, all the $g(i,n_i)$ that occur here are in $C$, and
therefore so are all the $\xi(i,n_i)$ (cf. \therosteritem2 and
\therosteritem4 in the
closure conditions for $C$) as well as $k_c$ (by closure condition
\therosteritem2).  Furthermore, by closure condition \therosteritem3,
$C$ also contains all the terms $g(i,n_i+1)$, so $y-k_y\in C\subseteq
D$.  But also $y\in D$, so we conclude $k_y\in D$.  Thus, the last
displayed equation above exhibits $d$ as 2 times an element $k_y$ of
$D\cap K$ plus terms from $C\cap K$.  This completes the proof that
$(D\cap K)/(C\cap K)$ is divisible and thus the proof of the theorem.
\qed\enddemo

\enddocument